\documentclass{article}
\ifx\pdfoutput\undefined % if we're not running pdftex
\else
\usepackage[pdftex,a4paper,colorlinks,hyperindex]{hyperref}
\fi
\usepackage{amsfonts}

\newtheorem{theorem}{Theorem}[section]
\newtheorem{lemma}[theorem]{Lemma}
\newtheorem{proposition}[theorem]{Proposition}
\newtheorem{corollary}[theorem]{Corollary}

\newenvironment{remark}{\refstepcounter{theorem}
 \bigbreak\noindent{\sc Remark }
\arabic{section}.\arabic{theorem}.}
 {\medbreak}

\newenvironment{definition}{\refstepcounter{theorem}
 \bigbreak\noindent{\sc Definition }
\arabic{section}.\arabic{theorem}.}
 {\medbreak}

\newenvironment{example}{\refstepcounter{theorem}
 \bigbreak\noindent{\sc Example }
\arabic{section}.\arabic{theorem}.}
 {\medbreak}

\newcounter{itemcounter}
\newenvironment{items}{
   \begin{list}{\alph{itemcounter})}
   {\usecounter{itemcounter}\setlength{\topsep}{3 pt}
   \setlength{\partopsep}{0 pt}\setlength{\itemsep}{0 pt}
      \setlength{\labelwidth}{2 em}
   }}{\end{list}}

\def\blockmatrix#1#2#3#4{\offinterlineskip
\left ( \vcenter{  \ialign{
\strut \ \hfil##\hfil \ &\vrule##\  &\ \hfil##\hfil \ \cr
\omit & height 2pt & \cr
$#1$ && $#2$ \cr
\omit & height 1pt & \cr
\noalign{\hrule}
\omit & height 2pt & \cr
$#3$ && $#4$ \cr }} \right)}

\def\cocoa
 {\mbox{\rm C\kern-.13em o\kern-.07em 
C\kern-.13em o\kern-.15em A}}

   \font\tengothic=eufm10
   \font\sevengothic=eufm7
   \font\fivegothic=eufm5

   \font\tenmsy=msbm10
   \font\sevenmsy=msbm7
   \font\fivemsy=msbm5

\newfam\msyfam
      \textfont\msyfam=\tenmsy
      \scriptfont\msyfam=\sevenmsy
      \scriptscriptfont\msyfam=\fivemsy
\newfam\gothicfam
   \textfont\gothicfam=\tengothic
   \scriptfont\gothicfam=\sevengothic
   \scriptscriptfont\gothicfam=\fivegothic

\def\goth#1{{\fam\gothicfam #1}}

\def\Cal#1{{\cal #1}}

\let\sem=\bf
\let\phi=\varphi
\let\rho=\varrho
\let\theta=\vartheta
\let\epsilon=\varepsilon

\def\TTo#1{\mathop{\longrightarrow}\limits ^{#1}}
\let\implies=\Rightarrow
\def\mapdown#1{\Big\downarrow\rlap{
$\vcenter{\hbox{$\scriptstyle #1$}}$}}

\def\squareforqed{\hbox{\rlap
{$\sqcap$}$\sqcup$}}
\def\qed{\ifmmode\squareforqed\else
{\unskip\nobreak\hfil
\penalty50\hskip1em\null\nobreak\hfil
\squareforqed
\parfillskip=0pt\finalhyphendemerits=0\endgraf}
\fi}

\def\proof{\rm \trivlist
 \item[\hskip \labelsep{\sc Proof.}]}
\def\endproof{\qed \endtrivlist}

\newcommand{\subsect}[1]{\refstepcounter{subsection}
\addcontentsline{toc}{subsection}
{\thesection.\Alph{subsection}\quad #1}
\subsection*{\thesection.\Alph{subsection}\quad #1}}

\def\LT{\mathop{\rm LT}\nolimits}
\def\init{\mathop{\rm in}\nolimits}
\def\gin{\mathop{\rm gin}\nolimits}
\def\GL{\mathop{\rm GL}\nolimits}
\def\g{\mathop{\goth g}\nolimits}

\def\matr{\mathop{\goth m}\nolimits}
\def\G{\mathop{\goth G}\nolimits}

\def\Ker{\mathop{\rm Ker}\nolimits}
\def\Im{\mathop{\rm Im}\nolimits}
\def\Mat{\mathop{\rm Mat}\nolimits}

\def\dprod{\mathop{\textstyle\prod}\nolimits}

\def\m{\mathop{\rm m}\nolimits}
\def\drl{{\mathop{\rm drl}\nolimits}}
\def\depth{\mathop{\rm depth}\nolimits}
\def\Sat{\mathop{\rm Sat}\nolimits}

\def\Stable{\mathop{\rm Stable}\nolimits}
\def\SStable{\mathop{\rm SStable}\nolimits}

\def\HF{\mathop{\rm HF}\nolimits}

\def\To{\longrightarrow}

\def\tsum{\mathop{\textstyle\sum}\limits}

\def\L{\mathop{{\mathcal{L}}}\nolimits}

\long\def \ignore#1{{}}

\title
{Generic Initial Ideals and Distractions
}

\author
{A.M.~Bigatti, A.~Conca, L.~Robbiano \thanks
{\tt AMS Classification: 13C99, 13P10}    
 }

\date{}

\begin{document}
\maketitle
%\label{firstpage}

\begin{abstract}
{\small
The generic initial ideals of a given ideal are rather recent
invariants. Not much is known about these objects, and it turns out to be 
very difficult to compute them. The main purpose of this paper
is to study the behaviour of generic initial ideals with respect to the  
operation of taking distractions. Theorem~\ref{gindl=I} is our  main result.
It states that  the {\tt DegRevLex}-generic initial ideal of the
distraction of a 
strongly  stable ideal is the ideal itself.  In proving this  fact we
develop some new results related to  distractions, stable and strongly
stable ideals. We draw some geometric conclusions for ideals of points.
}
\end{abstract}

%\tableofcontents

\bigskip\bigskip

\addcontentsline{toc}{section}{Introduction}

\section*{Introduction}

\begin{flushright}
\small\it 
What song the Syrens sang, or what name Achilles \\
assumed when he hid himself among women, \\
although puzzling questions are not\\
beyond {\sem all} conjecture.

\rm Sir Thomas Browne, Urn-Burial\\
(quoting the quote of E. A. Poe to \\
``The Murders in the Rue Morgue")
\end{flushright}

Everything related to the present paper started some time ago when 
we were playing with \cocoa. One should never underestimate the importance of
{\it playing}, even in Mathematics; fortunately nowadays we have
wonderful computer algebra {\it toys}. They are a sort of scientific
version of videogames  which for instance allow us to continually
explore the hidden secrets of polynomial rings. 

This paper provides another consequence of this new working behaviour
of many mathematicians. The story starts with  
the three of us  playing with generic initial ideals. First of all, it
should be said that when computing generic initial ideals, one's trust
in the answers  lies not only on the reliability of the computer
algebra system, but also on a more subtle belief. Namely, even for
very easy examples, it turns out that a {\it true} computation of
generic initial ideals is infeasible. When we say a true computation,
we mean a computation with a linear change of coordinates achieved with
a truly generic matrix, i.e.\ a matrix whose entries are distinct 
indeterminates.

What do we usually do instead? We perform a linear change of coordinates
with a {\it random}\/ matrix, and {\it trust}\/  the answer.
In practice, we believe that a random choice of a point in a
projective space does not hit a given hypersurface. Of course the
probability of hitting such a hypersurface is zero, but we might be
extremely unlucky. It is like thinking of a natural number and asking someone
else to guess it: the probability of success is zero, but ...
 
A~{\it good}\/ mathematician could argue that examples computed in this way 
are not acceptable, and from that point of view this paper
brings {\it good}\/ news, since we are able to produce a mathematical
proof which leads to the identification of an interesting class  
of generic initial ideals (see for instance Theorem~\ref{gindl=I},
Theorem~\ref{sumprinc}, and Corollary~\ref{ginandpoints}).

So, let us have a closer look at the content, and start off by 
recalling two known results about generic initial ideals. 

The first result (see Proposition~\ref{ginI=I})  says that if~$I$~is a
Borel-fixed ideal (see Definition~\ref{borelfixed}) in  the polynomial
ring
$P=K[x_1, \dots, x_n]$, and $\sigma$ is
a term ordering, then $\gin_{\sigma}(I) = I$.

The second result (see Theorem~\ref{ginh=ginxn}))  states that if~$I$~is
a homogeneous ideal in $P$, $h$ a  generic linear form,  and $\sigma$
a term ordering of $x_n$-DegRev type (see
Definition~\ref{defdegrevtype}), then we have 
$\gin_{\sigma_{\hat{x}_n}}(I_h) = \big(\gin_\sigma(I)\big)_{x_n}$.
Roughly speaking, the equality can be expressed by saying that if
$\sigma$ is  a term ordering of $x_n$-DegRev type, then taking the
$\gin_{\sigma_{\hat{x}_n}}$ of a generic hyperplane section of~$I$, is
like taking the $x_n$-hyperplane section of the $\gin_\sigma$ of~$I$.

While playing with \cocoa, we noticed that 
for a special kind of ideal, namely some ideals of distractions (see
Definitions~\ref{distrmatrix} and~\ref{distrpp})
stemming from strongly stable ideals (see Definition~\ref{ststable}), an
even more interesting equality holds, namely
$\gin_{\sigma_{\hat{x}_n}}(I_{x_n}) = \big(\gin_\sigma(I)\big)_{x_n}$, 
provided $x_n$ does not divide zero modulo~$I$. This fact was
purely empirical, but more than  enough to conjecture the
validity of the equality in the case of distractions of strongly
stable ideals and
$x_n$-DegRev type orderings. 

One might ask why we were considering strongly stable ideals. 
The motivation comes from the fact that these ideals are
Borel-fixed (see Proposition~\ref{charofborel}). And we have already
mentioned the fact that if $I$ is Borel-fixed, then $\gin_\sigma(I) =
I$ which might suggest other forms of stability for this type  of
ideal.
The turning point happened when
\cocoa\ completed the following session.

{\small
\begin{verbatim}
  Sigma    := Mat([[1,1,1,1], [0,0,0,-1], [1,0,0,0], [0,1,0,0]]);
  SigmaHat := Mat([[1,1,1],               [1,0,0],   [0,1,0]]);

  Use P ::= Q[x,y,z,w], Ord(Sigma);
  DI := ClassicDistraction(Ideal(x^5, x^4y, x^4z, x^3y^2, x^2y^3));

  Use PHat ::= Q[x,y,z], Ord(SigmaHat);
  DI_h := Image(DI, RMap(x, y, z, Randomized(x+y+z)));
  DI_w := Image(DI, RMap(x, y, z, 0));
  Gin(DI_h);   Gin(DI_w);
Ideal(x^5, x^4y, x^4z, x^3y^2, x^3yz, x^3z^3, x^2y^5)
-------------------------------
Ideal(x^5, x^4y, x^4z, x^3y^2, x^2y^3)
-------------------------------
\end{verbatim}
}

\noindent This was a counterexample to our conjecture.
Nevertheless we were unable to find any counterexample to the more
restricted conjecture  that the equality $\gin_{\sigma_{\hat{x}_n}}(I_{x_n}) =
\big(\gin_\sigma(I)\big)_{x_n} = I_{x_n}$ is true for the same kind
of ideals of distractions as above, and $\sigma ={\tt DegRevLex}$. 
A huge amount of
experimental evidence was accumulated. Fortunately, the solution,
although non trivial, turned out to be more reachable than discovering
what song the Syrens sang, and we were indeed able to prove the main
Theorem~\ref{mainth}, from which the validity of the conjecture stems as
an easy consequence, since Theorem~\ref{gindl=I} shows that 
$\gin_{_{\drl}}\big(D_{\L}(I) \big) = I$ for every strongly stable ideal~$I$.

At that point one might be lead to conjecture that 
$\gin_{_{\drl}}\big(D_{\L}(I) \big) = \gin_{\drl}(I)$
for every monomial ideal. This is not true even for stable ideals (see the 
second part of Example~\ref{counterex}). However, it is true for
special  monomial 
ideals which include {\it principal}\/ stable ideals (see
Theorem~\ref{sumprinc}). 

The road to the proofs had to be paved using some
preparatory material, and to this end we devote Sections 1, 2 and 3. In
particular, Section 2 develops some elements of the theory of
distractions. For instance, we give a
direct proof that monomial ideals and their distractions not only share
the same Hilbert function (Corollary~\ref{distrofideal}), but also share the
same standard Betti numbers (Theorem~\ref{preshomol} and
Corollary~\ref{distrandBettin}).

Distractions can be also described  as specializations of
polarizations of monomial ideals. They were (essentially)
introduced by  Hartshorne in his proof of the  connectedness of the
Hilbert scheme  \cite{H}, and since
then  they have been used, described, and rediscovered by many authors.
In this context let us mention the work of Pardue \cite{P1,P2} and the work of
Migliore and Nagel \cite{MN}.

\bigskip

Finally, we point out  an interesting application of our main
results. By  using suitable distractions of a strongly
stable monomial ideal 
$I$, one can  construct schemes of rational points whose 
$\gin$ with respect to ${\tt DegRevLex}$ is $I$ itself (see 
Corollary~\ref{ginandpoints}). In that case you do not have to rely on
luck to claim 
that you know the $\gin$, you have a proof!

%At this point one might be lured into trying to make a conjecture
%about the name 
%assumed by Achilles  when he hid himself among women. We have no clue, \cocoa\
%cannot help in that matter. 
%On the other hand, if you want
%to join us and play with the objects treated in  this paper, we
%recommend that you download \cocoa\ and start your own experimentation.

\newpage
\newpage

\section{Preliminaries}
\label{Preliminaries}

In this section we state the main preliminary facts. Although
all the results discussed here are known, the proofs given in 
the literature are not, in general, very accurate. For the sake  
of completeness we prefer to include some of them here. 

In the following
we let $P = K[x_1, \dots, x_n]$ be a polynomial ring with $n$
indeterminates over an {\sem infinite field} $K$.       We also
assume that it is standard graded, meaning that its grading 
satisfies
$\deg(c) = 0$ for every $c \in K\setminus \{0\}$, and 
$\deg(x_i) = 1$ for every 
$i = 1, \dots, n$. In other words, the degree matrix of the 
grading is $(1, 1, \dots, 1)$ (see \cite{KR2} for more details
on gradings  defined by matrices).

\medskip

The next definition introduces the notion of $\init_\sigma(f)$ 
for a  given term ordering~$\sigma$, and a non-zero polynomial
$f$. The full theory can be read for instance in~\cite{KR1}.
It~should be noted that here we use
$\init_\sigma(f)$ instead of
$\LT_\sigma(f)$, and the reason is that the notation
$\init(I)$ is compatible with $\gin(I)$, which will be the main 
object under investigation in this paper (gLT~does not look too
good, it does not sound too good either!).

\begin{definition}\label{initial ideal}
Let $I$ be an ideal in $P$ and let $\sigma$ be a term ordering
on $\mathbb{T}^n$, the monoid of power products in 
$x_1, \dots, x_n$. 
We denote by $\init_\sigma(I)$ the monomial ideal generated by 
$\{\init_\sigma(f) \mid f \in I\setminus\{0\} \}$ and call it
the  {\sem initial ideal}\/ of~$I$, or the {\sem leading term
ideal}\/ of~$I$.
\end{definition}

Throughout the paper we are always assuming that $\sigma$ is a
degree-compatible term ordering such that
$x_1 >_\sigma x_2 >_\sigma \cdots >_\sigma x_n$. The only
exception to this assumption is made when we talk about 
term orderings of $x_i$-DegRev type (see
Definition~\ref{defdegrevtype}, 
Proposition~\ref{propofxidegrev} and Theorem~\ref{ginh=ginxn}).

\medskip

Let us consider an automorphism ${\rm Aut}_K(P_1)$
on the $K$-vector space~$P_1$. If~we choose $x_1, \dots, x_n$ as 
a basis  of $P_1$, we may as usual describe the automorphism 
via a matrix 
$\g = (g_{ij}) \in \GL(n,K)$ with entries in $K$. 

If $\gamma \in  
{\rm Aut}_K(P_1) $ is represented by
$\g$, then $\gamma(x_j) = \sum_{i=1}^n g_{ij}x_i$. 
Therefore we have 
$$(\gamma(x_1), \dots, \gamma(x_n)) = 
(x_1,\dots, x_n)\cdot \g$$
where the ``$\cdot$" indicates the usual row-by-column product
of  matrices. The action of 
$\g$ on $P_1$ extends to $P$, so that it makes sense to consider
objects like
$\g(f)$ and $\g(I)$.

We can parametrize ${\rm Aut}_K(P_1)$  with
$\mathbb{A}_K^{n^2}\setminus H$, where $H$ is the
hypersurface defined by the vanishing of the determinant of  
a matrix whose entries are
$n^2$ distinct indeterminates. Having assumed $K$ to be
infinite, it make sense to speak of a {\it generic matrix},
likewise a {\it generic}\/ automorphism.

\medskip

\begin{definition}\label{borelfixed}
The subgroup $B$ of $\GL(n,K)$ of the upper triangular matrices
is called the {\sem Borel subgroup}.
If $I$ is an ideal in $P$, we say that $I$ is 
{\sem Borel-fixed} if  $\g(I) = I$ for every $\g \in B$.
\end{definition}

\begin{definition}\label{ststable}
A subvector space of~$P$ is called a {\sem monomial vector
space}\/ if it has a basis of power products.
An ideal in~$P$ which is a monomial vector space is called a
{\sem monomial ideal}.

If $t \in \mathbb{T}^n$, and $r$ is the biggest integer
such that $x_r\mid t$, we call $r$ the {\sem maximum index}\/ of
$t$ and denote it by $\m(t)$.
Let $T$ be a finite set of power products, and $V$ the
monomial vector space spanned by $T$. 
The set $T$ and the vector space
$V$ are called {\sem stable}\/ if $T$ satisfies the
following  property: for every power
product $t \in V$  and every $i\le \m(t)$, 
then $x_i \frac{t}{x_{\m(t)}} \in V$.  

The set $T$ and the 
vector space
$V$ are called {\sem strongly stable}\/ if $T$ satisfies the
following  property: for every power
product $t \in V$ such that $x_j\mid t$ and every $i\le j$, 
then $x_i \frac{t}{x_j} \in V$.  
\end{definition}

\begin{remark}\label{borelismon}
It is easy to see that if an ideal $I$ is Borel-fixed, then it 
is necessarily a monomial ideal.
\end{remark}

\begin{proposition}\label{charofborel}
Let $K$ be a field, and let~$I$ be a
monomial ideal in the polynomial ring $K[x_1, \dots, x_n]$.
Consider the following conditions.
\begin{items}
\item The ideal $I$ is Borel-fixed
\item The ideal $I$ is strongly stable
\end{items} 
Then $b) \implies a)$. Moreover, if ${\rm char}(K) = 0$ 
they are equivalent.
\end{proposition} 

\proof See for instance Proposition 1.25 in~\cite{G}.
\endproof

We are ready to recall a celebrated theorem. It was proved by
Galligo in characteristic zero, and then generalized by Bayer
and Stillman to every characteristic.

\begin{theorem}{\sem  (Galligo-Bayer-Stillman)}
\label{Galligo}\\
Let $K$ be an infinite  field,
let $\sigma$ be a term
ordering on
$\mathbb{T}^n$, and let $I$ be a homogeneous ideal in $P$.    
The following facts hold true.

\begin{items}
\item The ideal $\init_\sigma(\g(I))$ is constant for a generic $\g$.
\item The ideal $\init_\sigma(\g(I))$ is Borel-fixed for a generic $\g$.
\end{items}

\end{theorem}

\proof See for instance~\cite{G} p. 26. 
\endproof

This theorem allows us to give the following definition.

\begin{definition}\label{defofgin}
Let $\sigma$ be a term
ordering on
$\mathbb{T}^n$, and let $I$ be a homogeneous ideal in $P$. We 
denote by $\gin_\sigma(I)$ the Borel-fixed monomial ideal which
is equal  to $\init_\sigma(\g(I))$ for a generic $\g$, as
prescribed by the theorem. We call it the~{\sem generic
initial ideal}\/ of $I$ with respect to $\sigma$, or the {\sem
$\sigma$-generic initial ideal}\/ of $I$. Using this 
terminology, part b) of the theorem can be stated by saying  
that $\gin_\sigma(I)$ is Borel-fixed, and hence strongly   
stable in characteristic~$0$ by Proposition~\ref{charofborel}.
\end{definition}

\begin{proposition}{\sem (Gin and Borel-Fixed Ideals)}
\label{ginI=I}\\
Let $K$ be an infinite  field,
let $\sigma$ be a term
ordering on
$\mathbb{T}^n$, and let $I$ be a Borel-fixed ideal in $P$.    
Then $\gin_\sigma(I) = I$.
\end{proposition}
\proof It is a well-known fact, and here we give only a hint. 
Recall that  if
$\gamma
\in   {\rm Aut}_K(P_1)$ is represented by
$\g$, then 
$(\gamma(x_1), \dots, \gamma(x_n)) = 
(x_1,\dots, x_n)\cdot \g$.
If $\gamma$ is generic we can decompose $\g$ as the product of a lower
triangular 
times an upper triangular matrix.
Then we observe that the action of the lower triangular matrix
on a polynomial does not change its leading term. Moreover $I$
is closed under the action of the upper triangular matrix since~$I$
Borel-fixed.  
The
conclusion follows from the fact that $\gin_\sigma(I)$ and $I$ have the same
Hilbert function.
\endproof

\begin{example}\label{ginstabnonprinc}
The following is an easy example of a stable ideal which is not
strongly stable. 
In the polynomial ring $K[x_1,x_2,x_3]$, consider the 
following ideal 
$I=
(x_1^2,\ x_1x_2,\ x_2^2,\  \underline{x_2x_3})$. 
One can check (with \cocoa) that 
$\gin_{\drl}(I)=(x_1^2,\ x_1x_2,\ x_2^2,\ \underline{x_1x_3})$.
\end{example}

\medskip

At this point, we want to discuss what happens when we consider 
hyperplane sections.   Let $i \in\{1, \dots, n\}$, let $h
=\sum_{j=1}^n h_jx_j \in 
P_1$ be a linear form such that $h_i\ne 0$, and call 
$\hat P = K[x_1, \dots, x_{i-1},x_{i+1}, \dots, x_n]$.  
 The homomorphism $\phi: P \To \hat P$ defined by 
$\phi(x_j) = x_j$ for$j\ne i$, 
$\phi(x_i) = -\frac{1}{h_i}\big(\sum_{j\ne i}h_jx_j\big)$, 
induces an isomorphism  between $P/(h)$ and 
$\hat P$, which allows us to give the following definition.

\begin{definition}\label{hyperplaneconv}
Given a polynomial $f \in P$, we call $f_h$ the polynomial
$\phi(f)$ in
$\hat P$.
More generally, given a homogeneous ideal $I$ in $P$, we
call
$I_h$  the  ideal $\phi(I)$ in~$\hat P$. It is named
the {\sem $h$-hyperplane section ideal}\/ or simply the
\hbox{{\sem $h$-hyperplane section}\/ of $I$}.
\end{definition}

It should be observed that Definition~\ref{hyperplaneconv} is the correct
algebraic definition, but it does not correspond to the 
geometric notion of a hyperplane section. The reason is that, even if $I$ is
saturated,
$I_h$ need not to be such.

\medskip

Finally, we need another essential tool, namely the  notion of 
term  ordering of $x_i$-DegRev type. It is an important class of
term orderings which includes the so called {\tt
DegRevLex} ordering  (see~\cite{KR2}  for a through 
treatment of this topic). 

\begin{definition}\label{defdegrevtype}
Let $i\in\{1,\dots,n\}$. We say
that a  term ordering~$\sigma$ on $\mathbb{T}^n$ is of
{\sem {\boldmath $x_i$\unboldmath}-DegRev type},
if it satisfies the following conditions:

\begin{items}
\item The ordering~$\sigma$ is degree-compatible.
\item Given  $t,t' \in \mathbb{T}^n$, such that $\deg(t)=
\deg(t')$ and $\log_{x_i}(t) < \log_{x_i}(t')$,
we have $t >_\sigma t'$.
\end{items}

\end{definition}

We observe that term orderings of $x_i$-DegRev type exist. 
Namely, such an ordering can be obtained as {\tt Ord}$(V)$,
where $V$ is a matrix whose first two rows are 
{  $\tiny \pmatrix{1 &  \dots & 1 & \dots & 1 \cr
0 & \dots &-1 & \dots & 0\cr}$}.
         
\medskip

 The main result about term orderings of $x_i$-DegRev type is
contained in the following proposition. We use the notation 
$\sigma_{\hat{x}_i}$ to indicate the restriction of the term
ordering $\sigma$ to the monoid  of power 
products in the indeterminates 
$x_1, \dots, x_{i-1}, x_{i+1}, \dots, x_n$.

\begin{proposition}\label{propofxidegrev}
Let $I$ be a homogeneous ideal in $P$, let $i\in\{1,\dots,n\}$, 
let~$\sigma$ be a  term ordering of $x_i$-DegRev type
on~$\mathbb{T}^n$, and let $G=\{g_1,\dots,g_s\}$ be a
homogeneous
$\sigma$-Gr\"obner basis of~$I$.

\begin{items}
\item The set $ \{(g_1)_{x_i},\dots,  (g_s)_{x_i} \} \setminus
\{0\}$ is a homogeneous $\sigma_{\hat{x}_i}$-Gr\"obner basis of the
ideal~$I_{x_i}$.
\item We have $\big(\init_\sigma(I)\big)_{x_i} = 
\init_{\sigma_{\hat{x}_i}}(I_{x_i})$.
\item We have $\init_\sigma(I) +(x_i) = 
\init_{\sigma}\big(I +(x_i)\big)$.
\end{items}

\end{proposition}

\proof
The proof is easy. A full discussion and generalization can be 
found in~\cite{KR2}.
\endproof

Now we are ready to state an important theorem in the theory of 
generic initial ideals. 
Given a homogeneous ideal $I$, it compares the $\sigma_{\hat{x}_i}$-generic
initial ideal of the generic hyperplane section of~$I$, with the
$x_i$-hyperplane section of the $\sigma$-generic initial ideal
of~$I$. We were unable
to find any good proof   in the literature, so we decided to
include one here.

\begin{theorem}{\sem   (Gin and Hyperplane Sections)}
\label{ginh=ginxn}\\
Let $I$ be a homogeneous ideal in $P$, let $h \in P_1$ be a 
generic linear form, let $i\in\{1,\dots,n\}$, and let $\sigma$
be a term ordering of $x_i$-DegRev type. Then we have the
equality
$$\gin_{\sigma_{\hat{x}_i}}(I_h) = \big(\gin_\sigma(I)\big)_{x_i}$$
 of ideals in $K[x_1, \dots, x_{i-1},x_{i+1}, \dots, x_n]$.
\end{theorem}

\proof
To simplify the presentation, we give the proof for $i=n$,  since the 
general case is a straightforward generalization. 
The linear form $h =\sum_{i=1}^n h_ix_i$ is generic, hence 
$h_n \ne 0$.
We consider a generic matrix $\g \in \GL(n-1, K)$, and 
construct another matrix
{\tiny $\g' =  \blockmatrix \g {v^{\rm tr}} 0 0$},  where
$v = (-\frac{h_1}{h_n}, \cdots, -\frac{h_{n-1}}{h_n})$. 
We observe that  $v$ is  a generic vector in $K^{n-1}$.

It is clear that 
$$(x_1, \dots, x_{n-1},x_n)\cdot \g' = 
((x_1, \dots, x_{n-1})\cdot \g \ \mid \ 
-\frac{1}{h_n}\big(\tsum_{i=1}^{n-1}h_ix_i\big))$$
 hence we get 
the first equalities
$$\gin_{\sigma_{\hat{x}_n}}(I_h) =
\init_{\sigma_{\hat{x}_n}}\big(\g(I_h) \big) =  
\init_{\sigma_{\hat{x}_n}}(\g'(I))
\eqno{(1)}$$

Now, we consider a generic matrix $\G \in GL(n, K)$, and get 
the equality $\gin_\sigma(I) = \init_\sigma({\G}(I)) $ which
implies the equality
$$\big(\gin_\sigma(I)\big)_{x_n} = 
\big(\init_{\sigma}(\G(I))\big)_{x_n} \eqno{(2)}$$
We may  view $\G$ as a block matrix
in the following way
{\tiny  $\G = \blockmatrix \g {v^{\rm tr}} w c$}, where
$v, w$ are generic vectors in $K^{n-1}$,
and $c$ is a generic element in $K$.
We obtain
$$\big(\init_{\sigma}(\G(I))\big)_{x_n} =
\init_{\sigma_{\hat{x}_n}}\big(\G(I)_{x_n}\big) \eqno{(3)}$$ 
by Proposition~\ref{propofxidegrev}.b.
Using equalities (1), (2) and (3) it suffices to prove that
$\G(I)_{x_n} = \g'(I)$ which follows immediately from the 
definition of $\G$ and $\g'$. The proof is now complete.
\endproof

\begin{example}\label{ginxnneginxn}
With the following example we anticipate a theme which will be 
fully treated in the next sections.
As we said in the introduction, we asked ourselves the following
question. For ideals of distractions and a term ordering
$\sigma$ of $x_i$-DegRev type, is it true that 
$\gin_{\sigma_{\hat{x}_i}}(I_{x_i}) = \big(\gin_\sigma(I)\big)_{x_i}\!$?
\ In view of Theorem~\ref{ginh=ginxn}, the question could also be
phrased in the following way. Is it true that 
$\gin_{\sigma_{\hat{x}_i}}(I_{x_i}) = \gin_{\sigma_{\hat{x}_i}}(I_h)$?

There is a \cocoa\ counterexample, even if we start with a
strongly stable monomial ideal $I$. Namely, let $P= K[x, y, z,
w]$, and let $I$ be the ideal in $P$ generated by 
$\{x^5,\ x^4y,\ x^4z,\ x^3y^2,\ x^2y^3\}$ . 
Let $D(I)$ be the ideal in $P$ generated by
$\{ x(x-w)(x-2w)(x-3w)(x-4w), \ x(x-w)(x-2w)(x-3w)y,\ 
x(x-w)(x-2w)(x-3w)z,\ x(x-w)(x-2w)y(y-w),\ x(x-w)y(y-w)(y-2w)\}$.

Now let $\sigma = {\tt Ord}(W)$ be the
term ordering defined by the  matrix
$$W =   \pmatrix{1 & 1 & 1 & 1 \cr
0 & 0 & 0 & -1 \cr
1 & 0 & 0 & 0 \cr
0 & 1 & 0 & 0 \cr
}$$
 It is easy to see that $\sigma$ is of $w$-DegRev type, and
a 
\cocoa\ computation yields the following results.

\bigskip

\begin{tabular}{c c l}
$\gin_{\sigma_{\hat{w}}}(D(I)_h)$ & =  & $(x^5,\ x^4y,\ x^4z,\ x^3y^2,\
x^3yz,\ x^3z^3,\ x^2y^5)$ \\
\\
$\gin_{\sigma_{\hat{w}}}(D(I)_w)$  & =  & $\gin_{\sigma_{\hat{w}}}(I)$ \ = \
$I$ \\
\end{tabular}

\bigskip

\noindent Instead, if we do the same computation with $\sigma =
{\tt DegRevLex}$, we get that both $\gin_{\sigma_{\hat{w}}}(D(I)_h)$ and 
$\gin_{\sigma_{\hat{w}}}(D(I)_w)$ {\it coincide}\/ with $I$. We will see
later that this fact is {\it not a coincidence}.

\end{example}

WARNING: As we said in the introduction, 
\cocoa\ examples involving computations of generic initial
ideals, have probability of being correct equal to
$100\% - \epsilon$, where $\epsilon$ is as small as you wish...
but not equal 0.

\newpage

\section{Distractions}
\label{Distractions}

We are looking for special
situations where an equality of the type described in  
Theorem~\ref{ginh=ginxn} holds, but with $h$ replaced by
$x_n$. To achieve that goal, we need  more preparatory results.
In particular we need to build up some material around the
notion of distraction. To this end, we split the entire section
into three subsections, each addressing a specific question.

\subsect{First Properties}

We start out with a  definition which involves an
infinite  matrix of linear forms. The definition is given in
that way in order to simplify the discussion. It will be
clear, very soon, that in any application only a finite number
of columns play a role. As before, $K$ is an infinite field,
and $P = K[x_1, \dots, x_n]$.

\begin{definition}\label{distrmatrix}
We let $\L = \big (L_{ij} \mid 
i = 1, \dots, n,\ j \in \mathbb{N} \big)$
be an infinite matrix with entries $L_{ij} \in P_1$ with the
following two properties

\begin{items}
\item The equality $\langle L_{1\,j_1}, \dots,
L_{n\,j_n}\rangle = P_1$ holds for  every  $j_1, \dots,
j_n \in \mathbb{N}$. 
\item There exist an integer $N \in \mathbb{N}$ such that 
$L_{i\,j} =  L_{i\,N}$ for every $j > N$.
\end{items}  
 We call $\Cal{L}$ an {\sem $N$-distraction matrix}\/ or simply
a distraction matrix.
\end{definition}

\begin{definition}\label{distrpp}
Let $\L$ be a distraction matrix, and 
$t=x_1^{\alpha_1}x_2^{\alpha_2} \cdots x_n^{\alpha_n}$ a power 
product in $\mathbb{T}^n$.  Then the   polynomial 
$D_{\L}(t) = 
\dprod_{i=1}^n\big(\dprod_{j=1}^{\alpha_j}L_{ij} \big)$
is called the  {\sem $\L$-distraction}\/ of $t$.
Having defined $D_{\L}(t)$ for every power product $t$, we may
consider $D_{\L}$ as an operator acting on the power  products
which can then be extended by linearity. 
Therefore we can write $D_{\L}(V)$ where $V$ is a subvector
space of $P$, and call it the {\sem ${\L}$-distraction}\/ of
$V$.
\end{definition}

\begin{example}{\sem\  (Identical  Distractions)}
\label{identdistr}\\
 An obvious example of a distraction matrix is obtained by
taking $L_{ij} = x_i$ for every $i, j$. In this case the
operator
$D_{\L}$ is the identity operator. 
\end{example}

\begin{example}{\sem \ (Classic Distractions)}
\label{classicdistr}\\
 A classic example of distraction matrix is
given by choosing a big number
$N$, and then considering the assignment
$L_{ij} = x_i-(j-1)x_n$ for $i=1, \dots, n-1$ and  every 
$j < N$,  $L_{nj} = x_n$ for every $j \in \mathbb{N}$, and 
$L_{ij} = x_i$ for every $j \ge N$.
An important property of
this distraction is that $D_{\L}(t)$ is  square-free for every 
power product $t \in \mathbb{T}^{n-1}$, provided that $N$ and the
characteristic of $K$ are greater than the maximum exponent in $t$.
\end{example}

\begin{example}{\sem \ (1-Distractions)}
\label{genericdistr1}\\
It is interesting to observe that a $1$-distraction is
nothing else than a linear change of coordinates. Therefore,
if the linear forms $L_{11}, \dots, L_{n1}$ are generic, when
we apply the 1-distraction to a vector space we get a generic
linear change of coordinates.
\end{example}

\begin{example}{\sem \ (Generic Distractions)}
\label{genericdistr}\\
We observe that an infinite matrix with entries          
$L_{ij} \in P_1$ and property b) is a distraction matrix if for
every choice $j_1, \dots, j_n \in \mathbb{N}$, the square matrix
which represents $(L_{1\,j_1}, \dots, L_{n\,j_n})$ with respect
to the canonical basis $(x_1, \dots, x_n)$ is invertible.
We also observe that, due to property b), the number of
matrices required to be invertible is finite. Therefore, if 
the matrix $\L$ has property $b)$, and all its entries are
generic up to column $N$, then 
${\L}$ is an  $N$-distraction matrix.
\end{example}

Later on, we shall need this notion of genericity.

\begin{definition}\label{suffgeneric}
A distraction matrix ${\L}$ is said to be {\sem sufficiently
generic}\/ if $\langle L_{1\,j_1}, \dots,  L_{k\, j_k},
x_{k+1}, \dots, x_n \rangle = P_1$ for every $k =1, \dots, n$.
\end{definition}

A distraction matrix is sufficiently generic if 
for every choice $j_1, \dots, j_n \in \mathbb{N}$, all the
principal minors of the square
matrix which represents $(L_{1\,j_1}, \dots, L_{n\,j_n})$ with
respect to the canonical basis $(x_1, \dots, x_n)$, are 
invertible. Again, due to property b) of
Definition~\ref{distrmatrix}, the number of
minors required to be invertible is finite.

\medskip

We observe that  for
every power product
$t$, we have  $\deg(t) = \deg(D_{\L}(t))$, therefore, we can view
$D_{\L}$ as a
$K$-linear operator on every $P_d$. 
Moreover, if $t_1, t_2$ are power products such that
$t_1\mid t_2$, then 
$D_{\L}(t_1) \mid D_{\L}(t_2)$, but of course it is {\sem not
true} that  $D_{\L}(t_1t_2) = D_{\L}(t_1) D_{\L}(t_2)$, so
$D_{\L}:P \To P$ is not a homomorphism of rings.

\begin{proposition}\label{gLsuffgeneric}
Let ${\L}$ be a distraction
matrix, and let $\g \in \GL(n,K)$ be a matrix. 

\begin{items}
\item The matrix $\g\cdot {\L} =(\g\cdot \L_{ij})$ is a distraction
matrix.
\item  We have $\g D_{\L} = D_{\g\cdot {\L}}$.
\item If $\g$ is generic, then $\g \cdot {\L}$ is sufficiently generic.
\end{items}

\end{proposition}

\proof
Since a) is clear, we move to the 
proof of b).  It is sufficient  to prove the
equality 
$\g(D_{\L}(t)) = D_{\g \cdot {\L}}(t)$ for every power product $t$.
And such equality stems directly from the definitions.

Let us prove c). We have already observed right after the
definition of a sufficiently generic matrix that the number of
minors required to be invertible is finite. Let $\g$ be a matrix
whose entries are distinct indeterminates, and let 
$\matr_1, \dots, \matr_s$ be the square invertible
matrices which represent the $n$-tuples 
$(L_{1 j_1}, \dots, L_{n \,j_n})$ with respect to the canonical
basis $(x_1, \dots, x_n)$, and which are required to have
invertible principal minors. Now every matrix $\matr_i$ can
be put in triangular form with a suitable change of coordinates.
Therefore the condition on $\g$ to transform $\matr_i$ into a
matrix with the desired property is open and non-empty. Since
we have to achieve the property on a finite number of
matrices, we see that the condition on $\g$ is open and
non-empty. Now the proof is complete.
\endproof

The following proposition and subsequent corollary allow us to
extend the notion of 
${\L}$-distractions to ideals. If $f$ is a homogeneous
polynomial
$d$, we use the notation $fP_1$ to indicate the subvector space
$\langle fx_1, \dots, fx_n\rangle$ of $P_{d+1}$. More
generally, if $V$ is a subvector space of $P_d$ we use the
notation $VP_1$ to indicate the subvector space of $P_{d+1}$
spanned by the set
$\{vx_i \mid v \in V,\ i =1, \dots, n  \}$.

\begin{proposition}\label{distrtoideals}
Let ${\L}$ be a distraction matrix.  

\begin{items}
\item We have $D_{\L}(V\,P_d) = 
D_{\L}(V)\, P_d$ for every monomial 
vector space $V$, and every $d \in \mathbb{N}$.
\item The map $D_{\L}$ induces a $K$-linear automorphism of
$P_d$ for every $d \in \mathbb{N}$.
\item We have $\dim_K(V) = \dim_K\big(D_{\L}(V) \big)$  for
every subvector space $V$ of $P$.
\item  Let $V$, $W$ be subvector spaces of $P$. Then
$D_{\L}(V\cap W) = D_{\L}(V) \cap D_{\L}(W)$.
\item  Let $V$, $W$ be subvector spaces of $P\!\!$. Then
$D_{\L}(V\! +\! W) = D_{\L}(V)\! +\! D_{\L}(W)$.
\end{items}

\end{proposition}

\proof We
let $t=x_1^{\alpha_1}x_2^{\alpha_2} \cdots x_n^{\alpha_n}$ be a 
power product in
$\mathbb{T}^n$. Then 
$D_{\L}(t \,P_1) = \langle D_{\L}(tx_1), \dots,
D_{\L}(tx_n)\rangle =
\langle D_{\L}(x_1^{\alpha_1+1}x_2^{\alpha_2} \cdots 
x_n^{\alpha_n}), \dots, D_{\L}(x_1^{\alpha_1}x_2^{\alpha_2}
\cdots x_n^{\alpha_n+1}) \rangle = \langle D_{\L}(t)
L_{1,\alpha_1+1},
\dots,  D_{\L}(t) L_{n,\alpha_n+1}\rangle = D_{\L}(t)\, \langle 
L_{1,\alpha_1+1},
\dots, L_{n,\alpha_n+1} 
\rangle = D_{\L}(t)\, P_1$, where the last equality follows from
the assumption that ${\L}$ is a distraction matrix. We have
proved that $D_{\L}(t\,P_1) = D_{\L}(t)\, P_1$ which implies the
equality 
$D_{\L}(V\,P_1) = D_{\L}(V)\, P_1$. Now we prove the claim by
induction, namely we assume that $D_{\L}(V\,P_{d-1}) =
D_{\L}(V)\, P_{d-1}$, and get
$$D_{\L}(V\,P_d) = D_{\L}(V\,P_{d-1}\, P_1) =
D_{\L}(V\,P_{d-1})\,P_1 =  D_{\L}(V)\, P_{d-1}\, P_1 =
D_{\L}(V)\, P_d$$ The proof of  a) is now complete.

To prove b) we observe that, applying  a) to $V = K$, we get 
$D_{\L}(P_d) = P_d$, hence the $K$-linear operator 
$D_{\L}:P_d \To P_d$ is surjective and hence an automorphism.

The proofs of c), d) and e) follow directly from b).
\endproof

In the following, if $M$ is a standard graded $P$-module, we use
the notation
$\HF_M$ to indicate its  {\sem Hilbert function}.

\begin{corollary}{\sem  (Distractions and Hilbert Functions)}
\label{distrofideal}\\
Let ${\L}$ be a distraction matrix, and  $I$ 
a monomial ideal in $P$.

\begin{items}
\item  The vector space $D_{\L}(I)$ coincides with 
$\oplus_d \, D_{\L}(I_d)$, and is a homogeneous ideal in~$P$.
\item Let  $t_1, \dots, t_r \in \mathbb{T}^n$ be such that 
$I=(t_1, \dots, t_r)$. Then we have the equality $D_{\L}(I) =
\big(D_{\L}(t_1), \dots, D_{\L}(t_r)\big)$.
\item We have $\HF_I = \HF_{D_{\L}(I)}$.
\end{items}

\end{corollary}

\proof
The fact that $D_{\L}(I) = \oplus_d \, D_{\L}(I_d)$ follows from
the definition. To prove that it is a homogeneous ideal, we
have to show that
$D_{\L}(I_d)\,P_{d'} 
\subseteq D_{\L}(I_{d+d'})$ for every $d,d' \in \mathbb{N}$.
Indeed,
$D_{\L}(I_d)\,P_{d'} = D_{\L}(I_d\, P_{d'})$ by
Proposition~\ref{distrtoideals}.a, and the conclusion follows.
To prove b), we show that $D_{\L}(I)_d = \big(D_{\L}(t_1),
\dots, D_{\L}(t_r)\big)_d$ for every
$d \in \mathbb{N}$. To this end, we observe that $I_d = t_1\, 
P_{d-d_1} + \cdots +  t_r\, P_{d-d_r}$, where $d_j = \deg(t_j)$
for $j =1, \dots, r$. Therefore 
$$D_{\L}(I_d) = D_{\L}(t_1\, P_{d-d_1}) + \cdots D_{\L}(t_r\,
P_{d-d_r})= D_{\L}(t_1)\,P_{d-d_1} +\cdots +
D_{\L}(t_r)\,P_{d-d_r}$$ where the second equality follows from 
Proposition~\ref{distrtoideals}.a, and the proof is complete.
Finally, claim c) is an immediate consequence of
Proposition~\ref{distrtoideals}.c.
\endproof

The next result relates the operations of saturating and
distracting monomial ideals.
We use the notation $\Sat(I)$ to
indicate the {\sem saturation of} $I$, i.e.\ the unique ideal
$J$ such that
$\depth(P/J) > 0$ and $J_i = I_i$ for $i >>0$.

\begin{corollary}{\sem (Distractions and Saturation)}
\label{distrsat}\\
Let $\L$ be a distraction matrix and let $I$ be a monomial
ideal in $P$.  
Then $D_{\L}(\Sat(I)) = \Sat(D_{\L}(I))$.
\end{corollary}

\proof
We denote by $M$ the ideal $(x_1, \dots, x_n)$ and observe that 
$t \in I: M^r$ if and only if  $tP_r \subseteq I$ which is equivalent to 
$D_{\L}(tP_r) \subseteq D_{\L}(I)$ by Proposition~\ref{distrtoideals}.b.
Now $D_{\L}(tP_r) \subseteq D_{\L}(I)$ is equivalent to $D_{\L}(t)P_r
\subseteq D_{\L}(I)$  by Proposition~\ref{distrtoideals}.a. which is
equivalent to 
$D_{\L}(t) \in D_{\L}(I):M^r$, and the proof is complete.
\endproof

\subsect{Distractions and Betti numbers}

We have just proved that $I$ and $D_{\L}(I)$ share the same 
Hilbert function. Now we are going to prove a stronger
property, namely that 
$I$ and $D_{\L}(I)$ share also the same Betti numbers.

To achieve that result we need more facts about multigraded
rings and modules. A full treatment of this subject can be
found in~\cite{KR2}, from which we borrow some basic statements.

\begin{definition}\label{defgradedbyW}
Let $m\ge 1$, and let the polynomial ring $P=K[x_1,\dots,x_n]$
be equipped with a $\mathbb Z^m$-grading such that
$K \subseteq P_0$ and $x_1,\dots,x_n$ are homogeneous elements.

\begin{items}
\item For $j=1,\dots,n$, let $(w_{1j},\dots,w_{mj})\in  \mathbb
Z^m$ be the degree of~$x_j$. The matrix $W=(w_{ij})\in\Mat_{m,n}
(\mathbb Z)$ is called the {\sem degree matrix} of the given
grading. The rows of this matrix are called the {\sem weight
vectors} of the indeterminates $x_1,\dots,x_n$.

\item Conversely, given a matrix  
$W=(w_{ij})\in\Mat_{m,n}(\mathbb Z)$, we can consider the 
$\mathbb Z^m$-grading on~$P$ such that
$K \subseteq P_0$ and such that the indeterminates are
homogeneous elements whose degrees are given by the columns
of~$W$. In this case, we say that~$P$ is {\sem multigraded
by}~$W\!$, or simply {\sem graded by~$W$}.

\item Let $d\in\mathbb Z^m$.
The set of homogeneous polynomials of degree~$d$ is denoted
by $P_{W,d}$, or simply by $P_d$ if it is clear which grading
we are considering. A polynomial
$f\in P_{W,d}$ is also called homogeneous of {\sem
multidegree}~$d$, or simply of {\sem degree}~$d$,
and we write $\deg_W(f)=d$.

\end{items}
\end{definition}

\begin{proposition}\label{monismosthom}
Let~$I$ be an ideal of~$P$.
Then the following conditions are equivalent.

\begin{items}
\item The ideal~$I$ is a monomial ideal.

\item There is a non-singular matrix $W\in\Mat_n(\mathbb{Z})$
such that~$I$ is a homogeneous ideal
with respect to the grading on~$P$ given by~$W$.

\item For every $m\ge 1$ and every matrix
$W\in\Mat_{m,n}(\mathbb{Z})$, the ideal~$I$ is a homogeneous
ideal with respect to the grading on~$P$ given by~$W$.

\end{items}
\end{proposition}

\proof
%Since~$I$ is generated by terms and terms are homogeneous
%with respect to the gradings we are considering,
%Proposition~1.7.10 of~\cite{KR1} shows that~a) implies~c).
%Obviously,~b) is a special case of~c). Therefore is suffices to
%show that~b) implies~a).
%
%We take a homogeneous polynomial $f\in I_{W,d}$ and show that
%there is only one term in its support. Let $t=x_1^{\alpha_1}
%\cdots x_n^{\alpha_n}$ and $t'=x_1^{\beta_1}\cdots
%x_n^{\beta_n}$ be terms in the support of~$f$.
%Then $d=\deg_W(t)= \deg_W(t')$ implies
%$W\cdot (\alpha_1,\dots,\alpha_n)^{\rm tr}
%=W\cdot (\beta_1,\dots,\beta_n)^{\rm tr}$. Since
%we have $\det(W)\ne 0$, the $\mathbb{Z}$-linear map defined
%by~$W$ is injective. Hence  we obtain
%$(\alpha_1,\dots,\alpha_n)= (\beta_1,\dots,\beta_n)$, i.e.\ we
%get $t=t'$.

The easy proof is left to the reader. 
\endproof

This result allows us to say that a monomial
ideal $J$ is $\mathbb{I}$-homogeneous, where $\mathbb{I}$ is the identity
matrix. Therefore, we may construct a minimal finite free
resolution of $J$, which is also graded by $\mathbb{I}$. To understand
it,  it is therefore important to study
$\mathbb{I}$-graded modules and homomorphisms. 

Clearly, the $\mathbb{I}$-degree of a power product  $t$ in $P$
is
$\log(t)$, so that we can identify $t$ with its $\mathbb{I}$-degree,
and write $\deg_{\mathbb{I}}(t) = t$. Consequently, the graded free 
$P$-module generated by $1$ with $\deg_{\mathbb{I}}(1) = t$, can be
written as 
$P(-t)$, instead of  $P(-\log(t))$. 

 Next, we characterize 
$\mathbb{I}$-graded free modules, their $\mathbb{I}$-graded submodules,
and
$\mathbb{I}$-graded homomorphisms.

\begin{proposition}{\sem ($\mathbb{I}$-Graded Free Modules and
Homomorphisms)}
\label{charofzngrmod}\\
Let $P$ be graded by $\mathbb{I}$. 
\begin{items}
\item Let $\oplus_{i=1}^r P(-t_i)$ be a graded free
$P$-module, let $t \in \mathbb{T}^n$, and let $v$ be
a  homogeneous vector  of degree
$t$ in $\oplus_{i=1}^r P(-t_i)$. Then 
$v = (c_1s_1, \dots, c_rs_r)$ where $c_i \in K$, $s_i \in
\mathbb{T}^n$, and 
$s_it_i = t$ for every $i =1, \dots, r$.

\item Let $\oplus_{i=1}^r P(-t_i)$, 
$\oplus_{i=1}^{r'} P(-t_i')$ be finitely generated graded free
$P$-modules, and $\oplus_{i=1}^r P(-t_i) \TTo{\phi}
\oplus_{i=1}^{r'} P(-t_i')$ a graded $P$-homomorphism. If
$\phi$ is represented by a matrix~$\mathcal{M} = (m_{ij})$ with
respect to the canonical bases, then we have $m_{ij} = 0$ if
$t_i'$ does not divide $t_j$, and $m_{ij} = c_{ij}\, t_j/t_i'$
with $c_{ij} \in K$, if $t_i'$ divides $t_j$.

\end{items}
\end{proposition}

\proof
The proof of a) follows from the remark that  monomials
are the only homogeneous elements in $P$. To prove b), we observe that
$m_{ij}$ is the
$i^{\rm th}$ coordinate of $\phi(e_j)$. Now, $\deg(e_j) = t_j$,
hence $\deg(\phi(e_j)) = t_j$. But the only elements of degree
$t_j$ in $P(-t_i')$ are the monomials $c_{ij}\, t_j/t_i'$  if
$t_i'$ divides $t_j$. Instead, if $t_i'$ does not divide $t_j$,
only $0$ is homogeneous of degree $t_j$ in $P(-t_i')$.
\endproof

\begin{corollary}{\sem ($\mathbb{I}$-Graded Free Submodules)}
\label{idfreesubmod}\\
Let $P$ be graded by $\mathbb{I}$, let $\oplus_{i=1}^r P(-t_i)$
be an $\mathbb{I}$-graded free module, and let $M \subseteq
\oplus_{i=1}^r P(-t_i)$ be an
$\mathbb{I}$-graded submodule. Then $M$ is generated by a finite set of
homogeneous vectors, i.e.\ vectors  of the type  
$v = (c_1s_1, \dots, c_rs_r)$ where $c_i \in K$, 
$s_i \in \mathbb{T}^n$, and 
$s_1t_1 = s_2t_2 = \cdots = s_rt_r$.  
\end{corollary}

Let us look at Proposition~\ref{charofzngrmod}.c.
For the sake of simplicity,  we
are going to use the {\sem convention} that if $g$
divides $f$, and hence $f=gh$, the symbol $f\!::\!g$
means $h$,  otherwise it means~$0$. After the description given
in the above proposition, we allow ourselves to write 
\hbox{$\mathcal{M} = (m_{ij}) = (c_{ij}\, t_j\!::\!t_i')$.} 

\medskip

Let us turn to distractions.
We have seen that the operator $D_{\L}$ acts on
power products and extends to a $K$-linear map $D_{\L}:P \To P$
which is not a homomorphism of rings. 
To see how
it extends to $\mathbb{I}$-graded free $P$-mod\-ules, let
$\oplus_{i=1}^r P(-t_i)$ be an $\mathbb{I}$-graded free $P$-module, let
$d_i$ be the standard degree of $t_i$, and consider
the standard graded free
$P$-module $\oplus_{i=1}^r P(-d_i)$.

We have already observed that $t_1\mid t_2$ implies
$D_{\L}(t_1) \mid D_{\L}(t_2)$, and this fact allows us to give
the following definition.

\begin{definition}
We define a map from
$\oplus_{i=1}^r P(-t_i)$ to $\oplus_{i=1}^r P(-d_i)$
by sending every
power product $te_i$ to $D_{\L}(tt_i)/D_{\L}(t_i)\, e_i$. Then
we  extend it uniquely to a $K$-linear map  which
preserves the standard degrees, and which we  call
$D_{\L}$.  
\end{definition}

\begin{remark}
Arguing as in Corollary~\ref{distrofideal} and
its proof, we see that if $M$ is an $\mathbb{I}$-graded submodule
of $\oplus_{i=1}^r P(-t_i)$ (as described in
Corollary~\ref{idfreesubmod}), then 
$D_{\L}(M)$ is a standard graded submodule of $\oplus_{i=1}^r
P(-d_i)$. Moreover, they share the same Hilbert function.
\end{remark}

Next, we describe how the homomorphisms of $\mathbb{I}$-graded free
modules, and the maps~$D_{\L}$ work together.
\begin{lemma}\label{commsquare}
Let $\oplus_{i=1}^r P(-t_i)$ and 
$\oplus_{i=1}^{r'}P(-t_i')$ be
$\mathbb{I}$-graded free $P$-modules, let $d_i = \deg(t_i)$, $d'_i =
\deg(t_j)$ where $\deg$ is the standard degree, and let 
\hbox{$\phi: \oplus_{i=1}^r P(-t_i) \To
\oplus_{i=1}^{r'}P(-t_i')$} be a  homomorphism of $\mathbb{I}$-graded
$P$-modules. 
Assume that $\phi$ is represented by the matrix $\mathcal{M} =
(m_{ij}) = (c_{ij}\,t_j\!\!\!::\!\!\!t_i')$, let
$D_{\L}(\mathcal{M}) =  (c_{ij}\,D_{\L}(t_j)\!::\!D_{\L}(t_i'))$,
and
let $D_{\L}(\phi) : \oplus_{i=1}^r P(-d_i) \To 
\oplus_{i=1}^{r'} P(-d_i')$ be the standard graded
homomorphism of $P$-modules defined by 
$D_{\L}(\mathcal{M})$.
\begin{items}
\item The maps $D_{\L}$ are isomorphisms of $K$ vector spaces
which preserve the standard degree.
\item The maps $\phi$ and $D_{\L}(\phi)$ preserve the standard
degrees.
\item The diagram
$$\displaylines{ \hfill
\matrix{ 
\oplus_{i=1}^r P(-t_i) & \TTo{\phi} & 
\oplus_{i=1}^{r'}P(-t_i') \cr
\mapdown{D_{\L}} && \mapdown{D_{\L}} \cr
\oplus_{i=1}^r P(-d_i) & \TTo{D_{\L}(\phi)} &
\oplus_{i=1}^{r'} P(-d_i')  \cr     & & &}
\hfill}$$
is a commutative diagram of $K$-vector spaces.

\end{items}

\end{lemma}

\proof We have already observed that $D_{\L}$ preserves the 
standard degrees. It is an isomorphism by
Proposition~\ref{distrtoideals}.b. So the proof of a) is
complete.

We prove b).
The map $\phi$ is a
homomorphism of $\mathbb{I}$-graded $P$-modules, therefore it preserves
the fine degrees, and hence the standard degrees too. Finally,
we have $\deg \big(D_{\L}(t_j)|D_{\L}(t_i')\big) = d_j - d_i'$,
which implies that $D_{\L}(\phi)$ preserves the standard
degrees.

To show the commutativity of the diagram, it suffices to pick an
element of type $te_j \in \oplus_{i=1}^r P(-t_i)$ and show that 
$D_{\L}(\phi(te_j)) = D_{\L}(\phi)(D_{\L}(te_j))$, and to prove
such equality it suffices to prove it componentwise.
Let us compute the $i$-component of both vectors.
The $i$-component of $D_{\L}(\phi(te_j))$ is 
$$D_{\L}\big(c_{ij}t\, (t_j\!::\! t_i')\big) = 
c_{ij} D_{\L}\big(t_i'\;t\; (t_j\!::\! t_i')\big)/
D_{\L}(t_i')$$  and the latter is either 
$ c_{ij}\, D_{\L}(t\; t_j)/D_{\L}(t_i')$ or $0$. 
The $i$-component of the vector
$D_{\L}(\phi)\big(D_{\L}(te_j)\big)$ is the $i$-component of 
$D_{\L}(\phi) \big(c_{ij} D_{\L}(tt_j)/D_{\L}(t_j) e_j\big)$,
ie.\
$$
c_{ij}\,D_{\L}(tt_j)/D_{\L}(t_j)\; \big(D_{\L}(t_j)\!::\!
D_{\L}(t_i')\big)
$$ 
and the latter is either  
$ c_{ij}\,D_{\L}(tt_j)/ D_{\L}(t_i')$
 or $0$. The
proof is  complete.
\endproof

Using the above terminology, we are ready to prove an
interesting result.

\begin{theorem}{\sem (Distractions and
Homology)}\label{preshomol}\\
Let $\ \oplus_{i=1}^{r''} P(-t_i'')  \TTo{\psi} 
\oplus_{i=1}^{r}P(-t_i)   \TTo{\phi} 
\oplus_{i=1}^{r'}P(-t_i')\ $ be a complex of 
$\mathbb{I}$-graded free modules and homomorphisms.
Let $D_{\L}$, $D_{\L}(\psi)$,
$D_{\L}(\phi)$ be defined as above. 
\begin{items}
\item The diagram
$$\displaylines{ \hfill
\matrix{ 
\oplus_{i=1}^{r''} P(-t_i'') & \TTo{\psi} & 
\oplus_{i=1}^{r}P(-t_i)  & \TTo{\phi} &
\oplus_{i=1}^{r'}P(-t_i')\cr
\mapdown{D_{\L}} && \mapdown{D_{\L}} && \mapdown{D_{\L}} \cr
\oplus_{i=1}^{r''} P(-d_i'') & \TTo{D_{\L}(\psi)} &
\oplus_{i=1}^{r} P(-d_i) &\TTo{D_{\L}(\phi)} &
\oplus_{i=1}^{r'} P(-d_i')  \cr     & & & &}
\hfill}$$
is commutative.

\item We have $D_{\L}(\Ker(\phi)) = \Ker(D_{\L}(\phi))$, and 
 $D_{\L}(\Im(\phi)) = \Im(D_{\L}(\phi))$.
\item The map $D_{\L} : \oplus_{i=1}^{r}P(-t_i) \To
\oplus_{i=1}^r P(-d_i)$ induces a
$K$-linear map 
$$\overline{D}_{\L} : \Ker(\phi)/\Im(\psi) \To
\Ker(D_{\L}(\phi))/\Im(D_{\L}(\psi))$$
which is  compatible with the standard degree.
\item The map $\overline{D}_{\L}$ is an isomorphism of vector
spaces.
\end{items}
\end{theorem}

\proof Claim a) follows from 
Lemma~\ref{commsquare}.c.
The horizontal maps preserve the standard degrees by
Lemma~\ref{commsquare}.b., and the vertical maps are 
isomorphisms of vector spaces which preserve the standard
degree by Lemma~\ref{commsquare}.a.  The other conclusions
are achieved with standard arguments in homological algebra.
\endproof

Finally we can see the relationship between the Betti numbers
of $\mathbb{I}$-graded modules and their distractions.

\begin{corollary}{\sem (Distractions and Betti Numbers)}
\label{distrandBettin}\\
Let ${\L}$ be a distraction matrix, and  $M$ 
an $\mathbb{I}$-graded submodule of a finitely
generated $\mathbb{I}$-graded free module. The standard Betti numbers
of 
$M$ and $D_{\L}(M)$ are the same.
\end{corollary}

\proof
Let be given a minimal finite free $\mathbb{I}$-graded resolution of
$M$. We apply the operator $D_{\L}$ to the entire resolution,
and the theorem shows that we get a minimal finite free standard
graded resolution of $D_{\L}(M)$. Since the maps $D_{\L}$ are
isomorphisms of $K$-vector spaces, the standard Betti numbers
of $M$ and $D_{\L}(M)$ are indeed the same.
\endproof

%\begin{remark}
%{\tt TODO: completare \ } It would be possible to use the theory
%of polarization and show that a monomial ideal and a distraction
%are specialization of the same ideal. In this way, it would be
%possible to give another proof of the corollary above, at least
%for ideals.
%\end{remark}

\subsect{Radical Distractions}

The final part of this section addresses the following problem.
If $I$ is a monomial ideal, under which conditions on ${\L}$ is
$D_{\L}(I)$ a radical ideal?
 We recall that
{\sem irreducible monomial ideals} are of the type
$(x_{i_1}^{a_1}, \dots, x_{i_h}^{a_h})$, and that every
monomial ideal is the intersection of irreducible monomial
ideals, called {\sem irreducible components}.

\begin{definition}\label{radicdidtrmatr}
Let $I= (x_{i_1}^{a_1}, \dots, x_{i_h}^{a_h})$ be an
irreducible monomial ideal in $P$, and let 
\hbox{$S=\{(s_1, \dots, s_h) \mid 1\le s_1\!\le\!a_1, 
\dots, 1\le s_h \le a_h \}$}. 
Then  let
${\L}$ be a distraction matrix, and
let $V_{(s_1, \dots, s_h)}$ be the $K$-vector space
generated by
$\{  L_{i_1 s_1}, \dots, L_{i_h s_h} \}$.
If the vector spaces $V_{(s_1, \dots, s_h)}$ are pairwise
distinct when $(s_1, \dots, s_h)$ spans $S$, we say
that ${\L}$ is {\sem radical for $I$}. 
We observe that  a
necessary condition for an
$N$-distraction matrix to be radical for $I$ is that
$N \ge a_r$ for $r = 1, \dots, h$.

More generally, if $I$ is any
monomial ideal, we say that ${\L}$ is {\sem radical for} $I$ if
${\L}$ is radical for all the irreducible components of $I$.
\end{definition}

We observe that, with the exception of  the case $I =
(x_1,\dots, x_n)$, a necessary condition for a distraction
matrix to be radical for $I$ is that \hbox{ $\depth(P/I) > 0$},
since the vector spaces spanned by the linear forms taken one
from each of the rows coincide with $P_1$.

\begin{example} 
 The identical distractions of
Example~\ref{identdistr} are clearly not radical for any
irreducible monomial ideal, except the linear ones. 
Instead,  the
\hbox{$N$-generic} distractions of Example~\ref{genericdistr1}
are radical for every monomial ideal $I$~such that $\depth(P/I)
>0$, only provided that $N$ is sufficiently big.
\end{example}

\begin{example} 
Let us consider the  monomial ideal 
$I = (x_1^2x_2^2,\ x_1^2x_3^2,\ x_2^2x_3^2)$ in 
$P = K[x_1, x_2, x_3]$. Then
$$I = (x_1^2, x_2^2) \cap (x_1^2, x_3^2) \cap (x_2^2, x_3^2)
$$
is the decomposition of $I$ into irreducible components.
Now it is clear that the generic $2$-distraction is radical for
$I$, while the classic $2$-distraction is not.
\end{example}

Next result provides a reason for the word radical used in
the definition above.

\begin{lemma}\label{radirred}
Let $h<n$, let 
$I=(x_{i_1}^{a_1}, \dots, x_{i_h}^{a_h})$ be an irreducible
monomial ideal,  let
\hbox{$S=\{(s_1, \dots, s_h) \mid 1\le s_1\!\le
\!a_1, \dots, 1\le s_h \le a_h \}$}, and let ${\L}$ be a
distraction matrix which is radical for $I$. 
\begin{items}
\item We have $D_{\L}(I) = \bigcap_{(s_1, \dots, s_h) \in S}
(L_{i_1s_1}, \dots, L_{i_h s_h})$.
\item The ideal $D_{\L}(I)$ is radical.
\end{items}
\proof
Clearly b) follows immediately from a), and we prove a) by
induction on~$h$. The assumption on ${\L}$
implies that 
$L_{i_r 1}, L_{i_r2}, \dots, L_{i_ra_r}$ are pairwise coprime
for every $r =1, \dots, h$. Therefore if $h = 1$ the claim
follows. Let
$h>1$ and write
$I = (x_{i_1}^{a_1}) + J$, where $J = (x_2^{a_2}, \dots,
x_h^{a_h})$.
Then $D_{\L}(I) = (D_{\L}(x_1^{a_1})) + D_{\L}(J)$. By
induction, we have $D_{\L}(J) = \bigcap_{(s_2, \dots, s_h)}
(L_{2s_2}, \dots, L_{h s_h})$. Each ideal showing up in the
intersection is a linear, and hence prime, ideal of
height~$h-1$. Now, $L_{i_1 1}, \dots,L_{i_1 a_1}$ are pairwise
coprime and they do not divide zero modulo any of the prime
ideals in the intersection. Standard facts in ideal theory
imply that 
$D_{\L}(J) = \bigcap_{(s_1, \dots, s_h)}
(L_{1s_1}, \dots, L_{h s_h})$, and the proof is complete.
\endproof
\end{lemma}

\begin{remark}
An alternative way of proving the above lemma goes as follows.
Since 
$(x_{i_1}^{a_1}, \dots, x_{i_h}^{a_h})$ is a regular sequence,
a direct argument or an appropriate use of
Theorem~\ref{preshomol}, shows that $(D_{\L}(x_{i_1}^{a_1}),
\dots, D_{\L}(x_{i_h}^{a_h}))$ is a regular sequence too. This
fact implies that
$D_{\L}(I)$ is pure. Now all the ideals showing up in the
intersection are linear, hence prime of multiplicity $1$.
Clearly
$D_{\L}(I) \subset \bigcap_{(s_1, \dots, s_h)}
(L_{1s_1}, \dots, L_{h s_h})$, and both ideals are pure and
have the same dimension. From the assumption we deduce that
all the ideals in the intersection are pairwise distinct,
hence $D_{\L}(I)$ and $\bigcap_{(s_1, \dots, s_h)}
(L_{1s_1}, \dots, L_{h s_h})$ have the same multiplicity too. So
they are equal.
\end{remark}

\begin{proposition}{\sem (Distractions and Radical Ideals)}
\label{distisradical}\\
Let $I$ be a monomial ideal such that $\depth(P/I) >0$, 
and let ${\L}$
be a distraction matrix which is radical for $I$. Then the
ideal $D_{\L}(I)$ is radical.
\end{proposition}

\proof
We express $I$ as the intersection of irreducible ideals. The
assumption implies that every such irreducible ideal
$J$ is such that $\depth(P/J) >0$.
Therefore,
after Proposition~\ref{distrtoideals}.d  we may assume that
$I$ is irreducible and that $\depth(P/I) >0$, hence
that it is generated by
$h<n$ pure power products. At this point it suffices to apply
the lemma.
\endproof\newpage

\section{Stable and Strongly Stable Ideals}
\label{Stable and Strongly Stable Ideals}

In Section~\ref{Preliminaries} we introduced stable and strongly
stable ideals.  
Now we want to
investigate their structure more closely.

\begin{definition}
Let  $S=\{t_1,
\dots, t_s\}$ be a set of power products in $\mathbb{T}^n$. 
We denote by $\Stable(S)$ ($\SStable(S)$) the smallest
stable (strongly stable) ideal containing the elements of $S$. If
$S=\{t\}$ then we 
simply write $\Stable(t)$ for $\Stable(S)$, and $\SStable(t)$ for
$\SStable(S)$.  

\end{definition} 

\begin{definition}\label{elemmove}
If $t$ is a power product such that $x_j\mid t$ and $i\le j$,
we say that $x_i \frac{t}{x_j}$ is obtained from
$t$ with an {\sem elementary move}.  Likewise, we say that 
$x_i \frac{t}{x_{\m(t)}}$ is obtained from
$t$ with a {\sem special elementary move}.

Then we say that $t'$ is obtained
from $t$  with a {\sem chain of (special) elementary moves}, if there are
power products $t_1, \dots, t_s$ such that $t=t_1$, $t' = t_s$, and 
$t_i$ is obtained from $t_{i-1}$ with a (special) elementary move for $i
= 2, \dots, s$. 

Using this terminology, we see that for a (stable)
strongly stable set~$T$,  every power product reached from a
power product in $T$ with a chain of (special) elementary moves, still 
lies in $T$.
\end{definition}

It is  clear that if $I = \SStable(S)$ ($I = \Stable(S)$), every
minimal generator of $I$  
can be obtained by a chain of 
elementary moves (special elementary moves) from some elements in $S$. 
But we can say more.

\begin{lemma}\label{specstab}
Let $r \le n$, let $t_1,\dots, t_r \in \mathbb{T}^n$,  let  $\alpha_1,\dots,
\alpha_r$ be non-negative integers, and let $I$ be the ideal $\sum_{j=1}^r
t_j(x_1,\dots,x_j)^{\alpha_j}$. Assume that $t_1=1$,
$\m(t_j)<j$ for every $j=2,\dots, r$, that $t_j\, | \, t_{j+1}$  for
$1\leq j<r$, and that 
$\deg(t_j)+\alpha_j\leq \deg(t_{j+1})+\alpha_{j+1}$  for all
$1\leq j<r$. Then

\begin{items}
\item The ideal $I$ is  stable.
\item The graded Betti numbers and hence Hilbert
function of $I$ only depend on  the $r$-tuple of pairs
$((\deg(t_1),\alpha_1), \dots, (\deg(t_r),\alpha_r))$.
\end{items}
\end{lemma}

\begin{proof} Let $I_j=t_j(x_1,\dots,x_j)^{\alpha_j}$.  We claim that
$I$ is stable. For, let
$t$ be a power product in
$I$ and let
$j$ be the smallest integer such that
$t \in I_j$. Then $\m(t)=j$,   otherwise by
assumption we  would have $t \in I_{j-1}$. We
may write   $t=t_ju$, where $u$ is a power product  of degree
$\alpha_j$  such that $\m(u)=j$. Then clearly $t x_i/x_j$ is in $I_j$ and
hence in $I$  for all $i<j$.

A minimal free resolution of a stable monomial ideal $J$ is described by
Eliahou and Kervaire \cite{EK}. It follows from their work that the
graded Betti numbers of the  stable ideal $J$ are determined by the
multi-set of pairs of integers
$$
EK(J):=\{ (\m(t), \deg(t)) \ | \ t \mbox{ is a minimal generator of } J\}
$$
A minimal generator of $I$ is of the form $t_ju$ where $u$ is a
power product of degree $\alpha_j$, $\m(u)=j$, and  the
exponent of $x_j$ in $u$ is bigger than $\deg(t_j)+\alpha_j-
\deg(t_{j-1})-\alpha_{j-1}$. So the multi-set  $EK(I)$ only  depends on the set
of pairs   $(\deg (t_j),\alpha_j)$ with   $1\leq j \leq r$.
\end{proof}

\begin{proposition}\label{prinstab}
Let  $t=x_1^{a_1}\cdots x_n^{a_n} \in \mathbb{T}^n$ and let
$c_j=\sum_{k\geq j} a_k$ 
for $j=1,\dots, n$.     Then
$$\Stable(t)=\tsum_{j=1}^n x_1^{a_1}\cdots
x_{j-1}^{a_{j-1}}(x_1,\dots,x_j)^{c_j}$$
   \end{proposition}

\begin{proof} We let $I = \sum_{j=1}^n x_1^{a_1}\cdots
x_{j-1}^{a_{j-1}}(x_1,\dots,x_j)^{c_j }$. By taking 
$j=n$ we observe that $t$ belongs to  $I$.
Furthermore, it follows immediately from the definition of a stable ideal
that  $I$ is contained in $\Stable(t)$. The lemma shows that~$I$
is stable and the proof is complete.
\end{proof}

The next result yields an important characterization of strongly stable
ideals. Henceforth, we are going to call {\sem linear
segment ideals}\/ the ideals of the type $(x_1, x_2, \dots, x_h)$.

\begin{proposition} {\sem (Characterization of Strongly Stable Ideals)}
\label{charststid}\\
Let $I$ be a monomial ideal in $P$. The following conditions are
equivalent
\begin{items}
\item The ideal $I$ is strongly stable.
\item The ideal $I$ is a sum of ideals which are intersections of
powers of linear segment ideals.
\end{items}
\end{proposition} 

\proof
 It is clear that powers of linear
segment ideals are strongly stable and it is also clear that sums and
intersections of strongly stable ideals are strongly stable. 
On the other hand, let $I = (t_1, \dots, t_s)$ be a strongly stable
ideal. Clearly $I = \SStable(t_1) + \cdots + \SStable(t_s)$. To conclude the
proof, it suffices to observe that if $t = x_1^{\alpha_1}
x_2^{\alpha_2}
\cdots x_n^{\alpha_n}$ then $\SStable(t) = \cap_{i=1}^n (x_1, \dots,
x_i)^{\beta_i}$, where 
$\beta_i=\sum_{j=1}^i\alpha_j$.
\endproof

\begin{lemma}\label{satofstst}
Let $I$ be a strongly stable ideal such that $\depth(P/I) = 0$. 
Then $\Sat(I)$ is strongly stable.
\end{lemma}

\proof 
It is easy to see that $\Sat(I) = I:x_n^r$ for a sufficiently big $r$.
Let $t \in I$ be a power product, and let $t'$ be a power
product which is obtained from $t$ with a chain of elementary moves
(see Definition~\ref{elemmove}). Then $t x_n^r \in I$, and since
$t' x_n^r$ is also obtained from $t$ with a chain of elementary moves.
we get $t' x_n^r \in I$, hence $t' \in \Sat(I)$.
\endproof

\begin{example}
Consider the following
example. Let $I = (x,y^2)\cap (x^2, y, z)$ in $K[x,y,z]$. Then 
$I = (x^2, \ xy, \ y^2, \ xz)$, hence it is strongly stable. 
In agreement with the lemma
its saturation $(x, y^2)$ is strongly stable, but
its embedded component is not. However, it is also possible to write
$I = (x,y^2)\cap (x, y, z)^2$, and now both components are strongly
stable.
\end{example}

\newpage

\section{Gin and Distractions}
\label{Gin and Distractions}

This section contains the main results of the paper.
Before stating the main theorem, we fix a
bit of notation.  We  use the symbol $\init_{_{\drl}}$ to mean
$\init_{\tt DegRevLex}$. 
If $1\le m \le n$, we write $P_{[m]}$ instead of 
$K[x_1, \dots, x_{m-1}]$.
Moreover, if 
\hbox{${\L}=(L_{ij} \mid i=1, \dots, n, \ j \in \mathbb{N})$} 
is a distraction matrix, and $1\le m\le n$, we use
the notation 
${\L}_{[m]}$ to indicate the matrix $((L_{ij})_{[m]} \mid i=1,
\dots, m-1, \ j
\in \mathbb{N})$, where 
$(L_{ij})_{[m]}$ is the image  of 
$L_{ij} + (x_m, \dots, x_n)$ in $P_{[m]}$.

\medskip

It is clear that if ${\L}$ is a distraction matrix, then ${\L}_{[m]}$
need not be such. However, things change if ${\L}$ is a sufficiently
generic distraction matrix.

\begin{lemma} \label{suffgenrestr}
Let $\L = (L_{ij} \mid i=1, \dots, n, \ j \in \mathbb{N})$ 
be a sufficiently generic distraction matrix with entries in $P$, and
$1\le m\le n$.  Then ${\L}_{[m]}$ is
a sufficiently generic distraction matrix with
entries in
$P_{[m]}$.
\end{lemma}
\proof The claim follows directly from Definition~\ref{suffgeneric}.
\endproof

Now we are ready to prove the main result.

\begin{theorem} {\sem (Main Theorem)}
\label{mainth}\\
Let ${\L}$ be a sufficiently generic distraction matrix. Then 
$$\init_{_{\drl}}(D_{\L}(I)) = I$$
for every strongly stable monomial ideal $I$ in~$P$.
\end{theorem}

\proof
The proof is achieved through some intermediate claims.

\medskip

\noindent {\it Claim 1.} Let $I_1$, $I_2$ be two strongly stable
monomial ideals and suppose that the theorem is proved for $I_1$
and $I_2$. Then the theorem holds true for $I_1+I_2$.

\smallskip
{\it Proof of Claim 1}. We let $I = I_1+I_2$. From
Proposition~\ref{distrtoideals}.e we get the equality 
$D_{\L}(I) = D_{\L}(I_1) + D_{\L}(I_2)$.  We deduce the following chain
of relations
$$I = I_1 + I_2 =
\init_{_{\drl}}(D_{\L}(I_1)) + \init_{_{\drl}}(D_{\L}(I_2))  
 \subseteq 
\init_{_{\drl}}\big( D_{\L}(I_1)+D_{\L}(I_2) \big)  
$$
which shows that $I \subseteq \init_{_{{\drl}}}(D_{\L}(I))$. Since
the two ideals have the same Hilbert function by
Proposition~\ref{distrtoideals}.d, we may conclude that they
coincide, and Claim 1 is proved.

\medskip

\noindent {\it Claim 2.} Let $I_1$, $I_2$ be two strongly stable
monomial ideals and suppose that the theorem is proved for $I_1$
and $I_2$. Then the theorem holds true for $I_1\cap I_2$.

\smallskip
{\it Proof of Claim 2}. We let $I = I_1\cap I_2$. From
Corollary~\ref{distrofideal}.b we get the equality 
$D_{\L}(I) = D_{\L}(I_1) \cap  D_{\L}(I_2)$.  We deduce the following
chain of relations
$$I = I_1 \cap I_2  =  
\init_{_{\drl}}(D_{\L}(I_1)) \cap  \init_{_{\drl}}(D_{\L}(I_2))  
 \supseteq   \init_{_{\drl}}\big (D_{\L}(I_1) \cap D_{\L}(I_2) \big)
$$
We conclude again by invoking the coincidence of the Hilbert function
of the two ideals.

\medskip

\noindent {\it Claim 3.} Let $1\le m \le n$, let $I$ be a
monomial ideal in $P_{[m]}$, and assume that 
$\init_{_{\drl}}(D_{{\L}_{[m]}}(I)) = I$. Then 
$\init_{_{\drl}}\big(D_{{\L}_{[m+1]}}(IP_{[m+1]})\big) = IP_{[m+1]}$, where
$IP_{[m+1]}$ denotes the extension of $I$ to $P_{[m+1]}$.

\smallskip

{\it Proof of Claim 3}. Proposition~\ref{propofxidegrev}.c and
the obvious remark that {\tt DegRevLex} restricted to
$\mathbb{T}(x_1, \dots, x_m)$ is of
$x_m$-DegRev type imply the equality
$$\init_{_{\drl}}\big(D_{{\L}_{[m+1]}}(IP_{[m+1]})\big) +(x_m) = 
\init_{_{\drl}}\big(D_{{\L}_{[m+1]}}(IP_{[m+1]}) +(x_m)\big)
\eqno{(1)}$$ Then we observe that 
$$\init_{_{\drl}}\big(D_{{\L}_{[m+1]}}(IP_{[m+1]}) +(x_m)\big) =
\init_{_{\drl}}\big( (D_{{\L}_{[m]}}(I)) P_{[m+1]} +(x_m)\big)
\eqno{(2)}$$
since $I$ has all the generators in $P_{[m]}$. We use the same
argument as before, and say that 
$$\init_{_{\drl}}\big( (D_{{\L}_{[m+1]}}(I)) P_{[m+1]} +(x_m)\big)
=
\init_{_{\drl}}\big( (D_{{\L}_{[m+1]}}(I)) P_{[m+1]}\big) +(x_m)
\eqno{(3)}$$
and 
$$\init_{_{\drl}}\big( (D_{{\L}_{[m+1]}}(I)) P_{[m+1]}\big) +(x_m)
=  IP_{[m+1]} +(x_m) \eqno{(4)}$$
by assumption. Looking at the chain of equalities, we get 
$$\init_{_{\drl}}\big(D_{{\L}_{[m+1]}}(IP_{[m+1]})\big) +(x_m)  =
IP_{[m+1]} +(x_m)$$
and hence  $IP_{[m+1]} \subseteq
\init_{_{\drl}}\big(D_{{\L}_{[m+1]}}(IP_{[m+1]})\big)$. Since both
have the same Hilbert function by
Corollary~\ref{distrofideal}.c, it follows that they coincide,
and Claim 3 is proved.

\medskip

To finish the proof of the theorem, we observe that $I$ can be
expressed as a finite  combination of sums and intersections of powers
of linear segment ideals (see Proposition~\ref{charststid}).
Therefore, after Claims 1 and 2 it suffices to prove the theorem for
ideals of the type  $(x_1, \dots, x_m)^s$. Let us do it.
Consider the ring $P_{[m+1]}$,
and consider the ideal  $I = (x_1, \dots, x_m)^s$ in  $P_{[m+1]}$.
Such ideal coincides with $P_{[m+1]}$ from degree $s$ on. On the other
hand the assumption that ${\L}$ is a sufficiently generic
distraction matrix guarantees that its restriction ${\L}_{[m+1]}$ is
also a distraction matrix by Lemma~\ref{suffgenrestr}. So
$D_{{\L}_{[m+1]}}$ also coincides with $P_{[m+1]}$ from degree $s$ on,
by  Proposition~\ref{distrtoideals}.b. And then we get 
$$\init_{_{\drl}}\big(D_{{\L}_{[m+1]}}(I)\big) = I $$
We use Claim 3 to extend the above equality to the equality
$$\init_{_{\drl}}\big(D_{{\L}}(I P)\big) = I P $$
The proof is now complete.
\endproof

The caterpillar is ready to become a butterfly: we are now in a position
to tackle the conjecture described in the introduction and made so very
plausible by the enormous bulk of corroborative
computational evidence.

\begin{theorem}\label{gindl=I}
Let $I$ be a strongly stable monomial ideal in $P$, and let ${\L}$ be a
distraction matrix. Then 
$$\gin_{_{\drl}}\big(D_{\L}(I) \big) = I$$
\end{theorem}

\proof
We observe that $\gin_{_{\drl}}\big(D_{\L}(I)\big) =
\init_{_{\drl}}(\g(D_{\L}(I)))$ where $\g \in \GL(n,K)$ is generic.
Now we use Proposition~\ref{gLsuffgeneric}.b to get 
$\init_{_{\drl}}(\g(D_{\L}(I))) = \init_{_{\drl}}(D_{\g \cdot 
{\L}}(I))$, and  Proposition~\ref{gLsuffgeneric}.c to know that $\g \cdot 
{\L}$ is sufficiently generic. A direct application of the Main
Theorem~\ref{mainth} finishes the proof.
\endproof

The above formula might suggest that  $\gin_{_{\drl}}\big(D_{\L}(I)
\big) = \gin_{\drl}(I)$ 
for every monomial ideal. This is not true even for stable ideals (see the 
second part of Example~\ref{counterex}). However, it is true for
special monomial ideals 
which include {\it principal}\/ stable ideals, as we are going to show.

\begin{lemma}\label{ginspecstab}
Let $r \le n$, let $t_1,\dots, t_r \in \mathbb{T}^n$,  let  $\alpha_1,\dots,
\alpha_r$ be non-negative integers, and let $I$ be the ideal $\sum_{j=1}^r
t_j(x_1,\dots,x_j)^{\alpha_j}$. Assume that $t_1=1$,
$\m(t_j)<j$ for every $j=2,\dots, r$, that $t_j\, | \, t_{j+1}$  for
$1\leq j<r$, and that 
$\deg(t_j)+\alpha_j\leq \deg(t_{j+1})+\alpha_{j+1}$  for all
$1\leq j<r$. 

\begin{items}
\item We have $\gin_\sigma(I)=\tsum_{j=1}^r
  x_1^{\deg(t_j)}(x_1,\dots,x_j)^{\alpha_j}$ 
for every
term ordering $\sigma$ and for every base field $K$.
\item  We have 
$\init_{\drl}(D_{\L}(I))=  \gin_{\drl}(I)$ for every  sufficiently
  generic distraction~$\L$. 
\item We have $\gin_{\drl}(D_{\L}(I))=  \gin_{\drl}(I)$  for every
  distraction  $\L$  

\end{items}
\end{lemma}

\begin{proof} We let $\beta_j=\deg(t_j)$ for $j = 1, \dots, r$, 
$I_j=t_j(x_1,\dots,x_j)^{\alpha_j}$,
$J_j=x_1^{\beta_j} (x_1,\dots,x_j)^{\alpha_j}$, and $J=\sum_{j=1}^r
J_j$.

To prove a), we observe that  $J$ satisfies the assumptions of
Lemma~\ref{specstab}. Moreover, $J$ and $I$ share the same tuple
$((\deg(t_1),\alpha_1), \dots, (\deg(t_r),\alpha_r))$, and hence they
share the same Hilbert 
function by Lemma~\ref{specstab}.b. Knowing that, to prove the equality
$\gin_\sigma(I)=J$ it is enough to prove one inclusion, and we will show that
$\gin_\sigma(I)\supseteq J$.  Namely, from
$I=\sum_j I_j$ we get $\gin_\sigma(I)\supseteq \sum_j \gin_\sigma(I_j)$ and
$\gin_\sigma(I_j)=x_1^{\beta_j}\gin_\sigma((x_1,\dots, x_j)^{\alpha_j})=J_j$,
where the last equality follows from the fact that $(x_1,\dots,
x_j)^{\alpha_j}$ is strongly stable. The proof of~a) is now complete.

To prove b), we denote by $\L_j$ the distraction obtained from $\L$ by
shifting with 
$t_j$, i.e.\ $D_{\L_j}(t)=D_{\L}(t_jt)/D_{\L}(t_j)$. Then we have the equality 
$D_{\L}(I_j)=D_{\L}(m_j)D_{L_j}((x_1,\dots, x_j)^{\alpha_j})$, and hence
$$\init_{\drl}(D_{\L}(I))\supseteq \tsum_{j=1}^r \init_{\drl}(D_{\L}(I_j))=
\tsum_{j=1}^r x_1^{\beta_j}  \init_{\drl}( D_{\L_j}((x_1,\dots,
x_j)^{\alpha_j})$$
Moreover, we get $\init_{\drl}(
D_{\L_j}((x_1,\dots, x_j)^{\alpha_j})=(x_1,\dots, x_j)^{\alpha_j}$
as a consequence of Theorem~\ref{gindl=I}.
Summing up, we have shown that $\init_{\drl}(D_{\L}(I))\supseteq J$ and we
conclude that  $\init_{\drl}(D_{\L}(I))=J$ since the two ideals have
the same Hilbert function.

Claim c) follows from b) by Theorem~\ref{mainth}.
\end{proof}

We are ready to prove the promised result about {\it principal}\/
stable ideals.

\begin{theorem}\label{sumprinc}
Let $t=x_1^{a_1}\cdots x_n^{a_n} \in \mathbb{T}^n$ and $I=\Stable(t)$. Let
$b_j=\sum_{k<j} a_k$ for $j=1, \dots,n$, and $c_j=\sum_{k\geq j} a_k$.  
\begin{items}
\item We have
$\gin_\tau(I)=\sum_j  x_1^{b_j} (x_1,\dots,
x_j)^{c_j}$ for every
term ordering $\tau$ and every  base field $K$.
\item  We have $\gin_{\drl}(D_{\L}(I))=  \gin_{\drl}( I)$ for  every
  distraction  $\L$. 
\end{items}
\end{theorem}

\proof
It suffices to combine Proposition~\ref{prinstab} and  \ref{ginspecstab}.
\endproof

\begin{example}\label{counterex}
The following  examples  show that the above
theorem cannot be generalized to stable monomial ideals 
(see the warning given after Example~\ref{ginxnneginxn}).

In the polynomial ring $K[x_1,x_2,x_3,x_4]$, consider the 
following ideal 
$I=\Stable(\{x_1x_2,\ x_2x_3x_4\}) = 
(x_1^2,\ x_1x_2,\ x_2^3,\  x_2^2x_3,\ x_2x_3^2,\  x_2x_3x_4)$.
One can check (with \cocoa) that 
$$\gin_{\drl}(I)=(x_1^2,\ x_1x_2,\ x_2^3,\ x_2^2x_3,\ x_1x_3^2,\
x_2^2x_4)$$ 
$$\gin_{{\rm lex}}(I)=(x_1^2,\ x_1x_2,\ x_2^3,\ x_2^2x_3,\ x_1x_3^2,\
x_1x_3x_4)$$

\medskip

Let $I = \Stable(\{x_2^3,\ x_3^2x_4^2\})$, and $\L$ a generic
distraction. Then 
$$I = (x_1^3,\  x_1^2x_2,\    x_1x_2^2,\ x_2^3,$$ 
$$x_1^2x_3^2,\
x_1x_2x_3^2,\  x_1x_3^3,\  x_1x_3^2x_4,\ 
x_2^2x_3^2,\  x_2x_3^3,\  x_2x_3^2x_4,\  x_3^4,\  x_3^3x_4,\  x_3^2x_4^2)$$

$$\gin_\drl(I) =  (x_1^3,\  x_1^2x_2,\ x_1x_2^2,\   x_2^3,$$ 
$$x_1^2x_3^2,\  x_1^2x_3x_4,\  x_1^2x_4^2,\ 
  x_1x_2x_3^2,\  x_1x_2x_3x_4,\  x_1x_3^3,\ \underline{x_1x_3^2x_4},\ 
x_2^2x_3^2,\  x_2x_3^3,\  x_3^4)$$

$$\gin_\drl(D_{\L}(I)) = (x_1^3,\  x_1^2x_2,\ x_1x_2^2,\   x_2^3,$$ 
$$x_1^2x_3^2,\  x_1^2x_3x_4,\  x_1^2x_4^2,\ 
  x_1x_2x_3^2,\  x_1x_2x_3x_4,\  x_1x_3^3,\ \underline{x_2^2x_3x_4},\ 
x_2^2x_3^2,\  x_2x_3^3,\  x_3^4)$$

\end{example}

\bigskip

We observe that if the minimal generators of a strongly stable monomial
ideal~$I$ have a non-trivial
$\gcd$, then it has to be a pure power of $x_1$. So
we can generalize Main Theorem~\ref{mainth} a bit.

\begin{corollary}\label{ginwithgcd}
Let $I$ be a strongly stable monomial ideal in $P$, and let ${\L}$ be a
distraction matrix. Let $x_1^a$ be a non-trivial $\gcd$ of the minimal
generators of $I$, let $I = x_1^a\, J$, and let $F$ be a non-zero form
of degree~$a$.
\begin{items}
\item The monomial ideal $J$ is strongly stable.
\item We have $\gin_{_{\drl}}\big(F\, D_{\L}(J) \big) = I = x_1^a\, J$.
\end{items}
\end{corollary}

\proof
Claim a) follows directly from the definition of strongly stable ideal.
To prove b), we observe that after a generic change of coordinates,
$\init_{_{\drl}}(\g (F)) = x_1^a$. We argue as in the proof of
Theorem~\ref{gindl=I}, and get
$$\gin_{_{\drl}}\big(F\, D_{\L}(J)\big) =
\init_{_{\drl}}(\g (F\, D_{\L}(J)))= \init_{_{\drl}}(\g(F)\, 
(D_{\g \cdot  {\L}}(J))) 
$$
$$
= x_1^a\, \init_{_{\drl}}  (D_{\g \cdot  {\L}}(J)) = x_1^a J
$$
where the last equality follows directly from Theorem~\ref{gindl=I}
applied to~$J$.
\endproof

Next result compares the class of all
$\gin_{_{\drl}}(I)$  with the class of $\gin_{_{\drl}}(I)$
where $I$ is a radical ideal. 
%We call $\mathcal{G}in_{_{\drl}}$ the class of all the generic initial
%ideals with respect to {\tt DegRevLex}, and we call 
%$\mathcal{RG}in_{_{\drl}}$ the class of all the generic initial ideals
%of radical ideals
%with respect to {\tt DegRevLex}

\begin{corollary} {\sem (Gin and Radical Ideals)}
\label{ginandradical}\\
Let $I$ be a homogeneous ideal in $P$.
Assume that  $\depth(P/I) >0$ and that either ${\rm char}(K) = 0$
or ${\rm char}(K) > 0$ and 
$\gin_{_{\drl}}(I)$ is strongly stable. Then there exists a
homogeneous radical ideal $J$ in $P$ such that \hbox{$\,
\gin_{_{\drl}}(I) =
\gin_{_{\drl}}(J)$.}  
\end{corollary}

\proof
We let ${\L}$ be a generic $N$-distraction matrix with a sufficiently
big $N$. It is known that  $\depth(P/I) = \depth(P/\gin_\drl(I))$.
Proposition~\ref{distisradical} shows that
$D_{\L}(\gin_{_{\drl}}(I))$ is radical, and if $J =
D_{\L}(\gin_{_{\drl}}(I))$ we conclude the proof of Claim~a) by
Theorem~\ref{gindl=I}. 
Finally, Claim~b) is a consequence
of a) and Proposition~\ref{charofborel}.
\endproof

What happens if we drop the assumption that 
$\depth(P/I)>0$? Clearly, we may get around the difficulty by
embedding our ideal into a polynomial ring $\overline{P}$ with one more
indeterminate, and then apply Corollary~\ref{ginandradical} to $I$
considered as an ideal in $\overline{P}$.  The situation is slightly
more complicated if we do not want to extend $P$.

\begin{corollary} {\sem (Gin and Saturation)}
\label{ginandsat}\\
Let $I$ be a homogeneous ideal in $P$, and assume
that either ${\rm char}(K) = 0$
or ${\rm char}(K) > 0$ and 
$\gin_{_{\drl}}(I)$ is strongly stable.
\begin{items}
\item There exists a
homogeneous saturated  radical ideal $J$ in $P$ with the property that 
$\, \Sat\big(\gin_{_{\drl}}(I)\big) =
\gin_{_{\drl}}(J)$. 
\item There exists a
homogeneous saturated radical ideal $J$ in $P$ with the property that 
$\, {\rm Proj}\big(P/\gin_{_{\drl}}(I)\big) =
{\rm Proj}\big(P/\gin_{_{\drl}}(J)\big)$. 
\end{items}
\end{corollary}

\proof
Clearly Claim~b) is another way of formulating Claim~a). In order to
prove Claim~a) we distinguish two cases. If $\depth(P/I) > 0$, then
also 
$\depth(P/\gin_{_{\drl}}(I)) > 0$, and we conclude the proof by taking $J
=D_{\L}(\gin_{_{\drl}}(I))$ and  invoking
Corollary~\ref{ginandradical}.a.

Now let $\depth(P/I) = 0$, so that also $\depth(P/\gin_{_{\drl}}(I)) =
0$. In this case let $\gin_{_{\drl}}(I) = I' \cap Q$ be a decomposition
of the  ideal 
$\gin_{_{\drl}}(I)$ so that
$\depth(P/I') > 0$,
$\sqrt{Q} = (x_1,\dots, x_n)$, and hence 
$\Sat\big(\gin_{_{\drl}}(I)\big)= I'$.

 We let ${\L}$ be a generic $N$-distraction matrix with a sufficiently
big $N$.

\smallskip

\noindent {\it Claim 1.} The ideal $D_{\L}(I')$ is radical and is the
saturation of $D_{\L}\big(\gin_{_{\drl}}(I)\big)$.

\medskip
{\it Proof of Claim 1}. The fact that $D_{\L}(I')$ is radical 
follows from Proposition~\ref{distisradical}. The fact that it is the
saturation of $D_{\L}\big(\gin_{_{\drl}}(I)\big)$ follows from
Corollary~\ref{distrsat}.

\medskip
\noindent {\it Claim 2.} We have $I' = \gin_{_{\drl}}(D_{\L}(I'))$.

\smallskip
{\it Proof of Claim 2}. It suffices to use Lemma~\ref{satofstst} to
deduce that 
$I'$ is strongly stable, and then invoke Theorem~\ref{gindl=I}.

\medskip

\noindent We pick $J = D_{\L}(I')$, and observe that $J$ is radical by
Proposition~\ref{distisradical}  and saturated by Claim 1. Now 
$\Sat\big(\gin_{_{\drl}}(I)\big) = I'$, and $I' =
\gin_{_{\drl}}(J)$ by Claim~2, so that the proof is complete.
\endproof

As another consequence, we get the following corollary. This result brings  
finite sets of special points in  projective space into the
spotlight.

\begin{corollary} {\sem (Gin and Points)}
\label{ginandpoints}\\
Let $I$ be a zero-dimensional strongly stable monomial ideal in
$P$,  and let~${\L}$~be a distraction matrix which is radical for~$I$,
and whose entries are in the polynomial ring~$\overline{P}=K[x_1, \dots, x_n,
x_{n+1}]$.  Then 
$D_{\L}(I)$ is the ideal of a finite set of rational
points in $\mathbb{P}^n_K$ such that 
$\gin_{_{\drl}}\big(D_{\L}(I) \big) = I\overline{P}$. 
\end{corollary}

\proof
The equality $\gin_{_{\drl}}\big(D_{\L}(I) \big) = I\overline{P}$ follows 
directly from Theorem~\ref{gindl=I}.  So we have only to show that 
$D_{\L}(I)$ is the ideal of a finite set of rational
points in~$\mathbb{P}^n_K$. We have 
$$\dim(\overline{P}/D_{\L}(I)) = 
\dim \big(\overline{P}/\big(\gin_{_{\drl}}\big(D_{\L}(I) \big)\big) 
= \dim(\overline{P}/I) = \dim(P/I) + 1 = 1$$
The assumption on $I$ and Proposition~\ref{distisradical} imply
that $D_{\L}(I)$ is radical, hence
$D_{\L}(I)$ defines a radical zero-dimensional scheme in
$\mathbb{P}^n_K$, i.e.\ a scheme of points in 
$\mathbb{P}^n_K$.
They are rational, since their defining ideals are linear by
Lemma~\ref{radirred}.a, and the proof is complete.
\endproof

For example, if we are given a zero-dimensional strongly stable
monomial ideal $I$ in $K[x_1, \dots, x_n]$, we embed it in $K[x_1,
  \dots, x_n, x_{n+1}]$. 
Then we apply either the classic 
or the generic distraction with $N$ sufficiently large, and get an ideal
in $K[x_1, \dots, x_n, x_{n+1}]$ which defines a finite set of
rational points in 
$\mathbb{P}^n$, and whose
$\gin_{_{\drl}}$ is 
$I$ itself.
\newpage

\bigskip
Address of the authors: 

Dipartimento di Matematica,
Universit\`a di Genova, 

Via Dodecaneso 35, I-16146 Italy.

\bigskip

{\tt E-mail: bigatti,conca,robbiano@dima.unige.it~~~} 
\newpage

\end{document}